\documentclass{svmult}
\usepackage{color,graphicx,amssymb,amsmath,amsfonts,booktabs}

\usepackage{latexsym}
\usepackage{amsfonts}
\usepackage{amssymb}

\newenvironment{defin}{\begin{defn} \em}{\end{defn}}  

\newenvironment{theo}{\begin{theorem} \em}{\end{theorem}}
\newenvironment{prop}{\begin{propos} \em}{\end{propos}}

\newtheorem{defn}{{\normalsize\sc Definition}}[section]
\newtheorem{examp}{{\normalsize\sc Example}}[section]

\newtheorem{propos}{{\normalsize\sc Proposition}}[section]
\newtheorem{cor}{{\normalsize\sc Corollary}}[section]
\newtheorem{obser}{{\normalsize\sc Observation}}[section]
\newtheorem{alg}{{\normalsize\sc Algorithm}}[section]

\begin{document}

\title{Forcing Iterated Admissibility in Strategic Belief Models}
\author{
}
\institute{ 
\texttt{}}


\maketitle

\begin{abstract}
Iterated admissibility ({\em IA}) can be seen as exhibiting a minimal criterion of rationality in games. In order to make this intuition more precise, the epistemic characterization of this game-theoretic solution has been actively investigated in recent times: it has been shown that strategies surviving $m+1$ rounds of iterated admissibility may be identified as those that are obtained under a condition called {\it rationality and $m$ assumption of rationality} ($RmAR$) in {\em complete} lexicographic type structures. On the other hand, it has been shown that its limit condition, $R\infty AR$, might not be satisfied by any state in the epistemic structure, if the class of types is complete and the types are continuous. In this paper we introduce a weaker notion of completeness which is nonetheless sufficient to characterize IA in a highly general way as the class of strategies that indeed satisfy $R\infty AR$. The key methodological innovation involves defining a new notion of {\em generic types} and employing these in conjunction with Cohen's technique of {\em forcing}. 

\end{abstract}



\section{Introduction}
The epistemic analysis of game solutions has focused on the belief and knowledge requirements for their existence. Well-known solution concepts are evaluated in terms of how demanding those requirements are. One of the less demanding solutions is based on the idea that rational players will never choose {\em weakly dominated} ({\em inadmissible}) actions \cite{PR}. Furthermore, the ideal of a {\em strategically stable} equilibrium can be achieved through the iterated elimination of inadmissible strategies \cite{KM}. Profiles that capture this intuition are composed by actions that survive the iterated deletion of weakly dominated strategies at successive levels. Strategies that meet this condition are known as {\em iterated admissible} ({\em IA}) strategies.  

The epistemic characterization of {\em IA} has been actively investigated in recent times, because of the challenges that it poses. On one hand, the elimination of a weakly dominated strategy at a given round of deletion requires the assumption that strategies deleted in previous rounds will not be used. Yet on the other, the definition of admissibility indicates that at whatever relevant level of iteration some probability distribution exists in the context of which a certain admissible strategy registers as a best response relative to all other possibilities. The conceptual support for this notion thus entails the potential inclusion of all the strategies of the other players. Therefore, none can be discarded. 

This {\em inclusion-exclusion} problem has been formally disposed by resorting to {\em Lexicographic Probability Systems} \cite{BFK}. A key result in this approach is that strategies that survive $m+1$ rounds of iterated admissibility can be obtained if players satisfy the epistemic condition known as {\it rationality and $m$ assumption of rationality} ($RmAR$) in {\em complete} type structures \cite{BFK}. In other words, intuitively, if there are ``enough'' types in a given model then a movement preserving admissible solutions from level $m$ to level $m+1$ is always possible in that model under the $RmAR$ condition. This would seem to suggest that the limit case at infinity would also provide a solution (in a precise sense that we clarify below), but it has been shown that this is not the case, that is, $R\infty AR$ might not be satisfied by any state in the epistemic structure, if the types are continuous and the class of types is complete. Nonetheless, if continuity is not required \cite{KL} or completeness is not demanded {\cite{P} it can be proven that a profile of strategies does indeed satisfy $R\infty AR$ if and only if it satisfies $IA$.

In this paper we investigate issues corresponding to the above in {\em Qualitative Type Spaces} \cite{PR}. Within this latter framework, it has in fact been shown that the construction of a complete type structure \cite{BK} is impossible. Our strategy is -- despite the impossibility of constructing a complete type structure -- to make use in this framework of the natural order relations emerging from the hierarchy of iterated levels in order to generate a ``new kind'' of type, called {\it generic}, which, when adjoined to the initial model, will be sufficient to broaden it and thereby to guarantee that $R\infty AR$ is satisfied within it by a straightforward application of Cohen's technique of {\em forcing}. This, in conjunction with the result indicated above, will in turn entail that $IA$ is also satisfied in the enlarged (although not complete) model.

The mathematical tool of forcing was originally introduced by Paul Cohen in order to show that Georg Cantor's famous Continuum Hypothesis is independent of the axioms of Zermelo-Frenkel Set Theory \cite{Co} \cite{K}. Since then, forcing has been mostly relegated to the realm of the foundations of Mathematics; its use remains uncommon in applied fields. Philosophically, this may be understood as a consequence of Shoenfield's Theorem, which entails that forcing yields results only in the non-absolute fragment of Mathematics, while most of applied science would seem to be confined in the absolute realm \cite{Je}. In this work, however, we show that forcing can indeed be used in the present case to provide definite features to the types associated to $IA$ strategy profiles. 

This paper is organized as follows:  In section 2 we introduce the notion of $IA$ solution while in section 3 we present a possibility-based epistemic model of games and strategies. In section 4 a brief conceptual discussion of Cohen's method of forcing is given. In section 5 we then use forcing to address the problem of defining {\em generic types} in strategic belief models. Finally, in section 6 we show how these types ensure that $R\infty AR$ is equivalent (in terms of strategies) to $IA$ in this new context.

\section{Iterated Admissibility}

The strategic interaction among self-interested parties can be represented as a {\em strategic form game}:

\begin{defin}
Let $G = \langle I, \{S_i\}_{i \in I}, \{\pi_i\}_{i \in I} \rangle$ be a {\em game}, where $I=\{1,\ldots, n\}$ is a set of {\em players} and $S_i, i\in I$ is a finite set of {\em strategies} for each player. A {\em profile} of strategies, $s=(s_1, \ldots, s_n)$ is an element of $S = \prod_{i \in I} S_i$.
In turn, $\Pi_i : S \rightarrow R $ is player $i$'s payoff.
\end{defin}

The goal is to assess the {\em solutions} of the game, i.e. the family of profiles $\mathcal{S} \subseteq S$ that might be expected to be chosen by the players. These profiles capture the rationality of players, seeking to maximize their payoffs. Being $\Pi_i(s)$ dependent on the choice of $i$, $s_i$ as well as of the choices of the rest of the players (denoted $s_{-i} \in \prod_{j \neq i} S_j$), a profile $(s_i, s_{-i}) \in \mathcal{S}$ shows also the coordination among agents that ensures this outcome. 

One of the many possible solution notions stems from a basic postulate of decision theory, namely that no player will choose a weakly dominated ({\em inadmissible}) strategy. That is, if for $s_i \in S_i$ there exists $\sigma^{*}_i \in \Delta S_i$ (i.e. a probability distribution over $S_i$ such that $\Pi_i(\sigma^{*}_i, s_{-i}) \geq \Pi_i(s_i, s_{-i})$ for every $s_{-i} \in S_{-i}$ and $\Pi_i(\sigma^{*}_i, \bar{s}_{-i}) > \Pi_i(s_i, \bar{s}_{-i})$ for some $\bar{s}_{-i} \in S_{-i}$, $s_i$ will not be chosen. Otherwise, $s_i$ is called {\em admissible} with respect to $S_i \times S_{-i}$. Then, if we set $S_i^{0} = S_i$ for every $i$, we can define inductively:

$$S_i^{m+1}\ = \ \{s_i \in S_i^{m}: s_i \ \mbox{is admissible with respect to} \ \prod_{i=1}^{n}S_i^{m} \}$$   

\noindent where $s_i \in S_i^{m}$ is said $m$-admissible, while $s_i \in S_i^{\infty} = \bigcap_{m=0}^{\infty} S_i^{m}$ is {\it iteratively admissible}. Thus, the set of {\em iterated admissible} ($IA$) strategies yields the solution notion $\mathcal{S}_{IA} = \prod_i^n S_i^{\infty}$.

The interest in $\mathcal{S}_{IA}$ resides in its naturality and in the achievement of {\em strategic stability} it ensures \cite{KM}. In turn, it raises the question of the aspects that contribute to the coordination on $\mathcal{S}_{IA}$, which is implicit since no communication is allowed among the players. In particular, we might inquire into those aspects of an {\em epistemic} nature, that is, those ensuing from the beliefs and knowledge of the individual players.  

\section{Strategic Belief Model of a Game}

We now show how these notions apply within {\em Qualitative Type Spaces}, which allow us to characterize certain strategic aspects of games.  In particular, we define the notion of {\em Strategic Belief Model}, which captures  the epistemic aspects involved in any choice of strategies \cite{BK}:

\[{\mathbf M}= (\{S_i\}_{i\in I}, \{T_i\}_{i\in I}, \{P_i\}_{i\in I})\]

\noindent where for each $i \in I$, $S_i$ and $T_i$ are $i$'s sets of strategies  and {\em types}, respectively. The structure intends to model a  game $G$ and each strategy-type pair is a {\em state} for a player, and each type of a player has beliefs about the states of the other players. These beliefs are captured by the relations $P_i$ that satisfy:
\begin{enumerate}
  \item $P_i: T_i \rightarrow S_{-i} \times T_{-i}$ is a correspondence.
  \item For all $t_i \in T_i$, $P_i[t_i] \neq \emptyset$.
\end{enumerate}

That is, $P_i[t_i]$ captures the strategies and types of the other players that $i$ thinks are possible. 

The analysis of the rationality of players and the epistemic conditions of solutions to the game can be fully disclosed up from the {\em states of the game}, i.e. profiles of states of the players. 

$\mathbf{M}$ has interesting features. The first one is that each $t_i$ can be ``unfolded'' in terms of the types of the other players, which in turn leads to beliefs about the type of $i$, etc. 

To make this notion more precise, let us define for any $t_i$ of $i$, the unfolding of $t_i$:

\begin{itemize}
\item $t_j \in P^{1}_{i}[t_i]$ if there exists $s_j \in S_j$ with\footnote{Here $s_{-ij}$ ($t_{-ij}$) denotes an element in $\prod_{k \neq i, k \neq j} S_k$ ($\prod_{k \neq i, k \neq j} T_k$).} $\langle (s_{-ij}, t_{-ij}),(s_j, t_j) \rangle \in P_i[t_i]$.
\item $t_j \in P^{m}_i[t_i]$, for any natural number $m$, if there exists $t_k \in P^{(m-1)}_i[t_i]$ such that $t_j \in P^{1}_k[t_k]$. 
\end{itemize}

This means that, if $t_j \in P^{m}_{i}[t_i]$, $t_j$ can be unfolded in $m$ steps to $t_i$, i.e. $t_j$ is believed by $t_i$ by considering $m$ steps of belief.

Another important feature of $\mathbf{M}$ is that it provides a powerful framework for describing the epistemic aspects involved in a game \cite{AB95}. The fundamental concept here is that of {\em assumption} defined over {\em events} of $\mathbf{M}$, i.e. on sets of states of the game. For any $E \subseteq \prod_{i} (S_i \times T_i)$, the types of $i$ that assume $E$ are denoted as ${\mathcal A}{\mathcal S}_{i}[E]$ with:\footnote{The notation $E_{|\prod_{j\neq i} (S_j \times T_j)}\}$ indicates the projection of $E$, defined over $\prod_{i} (S_i \times T_i)$ over $\prod_{j\neq i} (S_j \times T_j)$.}

$${\mathcal A}{\mathcal S}_{i}[E]= \{ t_i \in T_i: P_i[t_i] = E_{|\prod_{j\neq i} (S_j \times T_j)}\}$$

In particular, we say that at $t_i$, $i$ assumes that $j \neq i$ is rational if $t_i \in {\mathcal A}{\mathcal S}_{i}[j \ \mbox{is rational}]$, where the event ``$j$ is rational'' is
$$\{ \langle s_j, t_j; s_{-j}, t_{-j}\rangle: \Pi_j (s_j, s_{-j}) \geq \Pi_j (s, s_{-j}) \ \mbox{for any} \ s \in S_j \ \mbox{with} \ s_{-j} \in P_j[t_j]_{|S_{-j}}\}$$ 

Then, we can define the condition denoted $RmAR$ (for {\em Rationality and} $m${\em -Assumption of Rationality}):\footnote{If $E^i \subseteq S_i \times T_i$, ${\mathcal A}{\mathcal S}_{i}[E^i]$ is a shorthand for ${\mathcal A}{\mathcal S}_{i}[\bigcap_{E: E_{|S_i \times T_i} = E^i} E]$.}

\begin{itemize}
\item $R^{0}_{i}= [i \ \mbox{is rational}]_{|S_i \times T_i}$.

\item $R^{m}_{i}= R^{m-1}_{i} \cap (S_i \times {\mathcal A}{\mathcal S}_{i}[\bigcap_{j \neq i} R^{m-1}_{j}])$.
\end{itemize}

Then, finally, we can define condition $R\infty AR$ as that satisfied by all the states in  $\prod_{i}^n R^{\infty}_i$, where $R^{\infty}_i = \bigcap_{m \geq 0} R^{m}_i$ for every $i$. That is, we have formally characterized all the states in which players are rational, assume the rationality of the others, assume that the others assume their own rationality, etc.

\section{$R\infty AR$ and $IA$}

It is quite natural to think that $R\infty AR$ is an epistemic condition that supports $\mathcal{S}_{IA}$. Indeed, this relation can be strengthened to make them equivalent, but this identification hinges on the properties of the space of types. More particularly it relies on its richness, since there have to be ``enough'' types to complement the strategies that survived the IA process and make them states in $R\infty AR$. 

Among the conditions that ensure the richness of the type space is {\em completeness}, i.e. that for each player $i$ and every $E \subseteq S_{-i} \times T_{-i}$ there exists a type of $i$, say $t_i$ such that $P_{i}[t_i]= E$. But as shown in \cite{BK}, a strategic belief model cannot in fact satisfy this condition. That is, there necessarily exist events that cannot be assumed by the players.

A weaker condition but one still strong enough may be provided as follows: given a strategy $s_i^* \in S_i$, admissible with respect to $\bar{S}_i \times \bar{S}_{-i} \subseteq S_i \times S_{-i}$ there exists some $t_i \in T_i$ such that $\Pi_i(s^*_i, s_{-i}) \geq \Pi_{i}(s_i, s_{-i})$ for any $s_{-i} \in P_{i}[t_i]_{|\bar{S}_{-i}}$ and any $s_i \in \bar{S}_i$. If this condition is satisfied by $\mathbf{M}$ for any $i$ and every $\bar{S}_i \times \bar{S}_{-i} \subseteq S_i \times S_{-i}$, $\mathbf{M}$ may be said to be {\em rationality-complete}.

We then have the following result:

\begin{theo}\label{admiss}
If $\mathbf{M}$ is rationality-complete and $\prod_{i}^n R^{\infty}_i$$\neq$$\emptyset$, then $R^{\infty}_{i|S_i} = S_i^{\infty}$ for every $i$.
\end{theo}
\noindent {\bf Proof}: {\it We will show that for any $i$, if $R^{\infty}_i \neq \emptyset$ then for every $m$, $R^{m}_{i|S_i} = S_i^{m+1}$. We will proceed by induction over $m$}:
\begin{itemize}
\item $m=0$) {\it Consider $s_i \in R^{0}_{i| S_i}$. Then, there exists $t_i$ such that $(s_i,t_i)$ is rational. Suppose $s_i$ is not admissible with respect to $S_i \times S_{-i}$. Then there exists $\hat{s}_i \in S_i$ such that $\Pi_i (\hat{s}_i , s_{-i}) \geq \Pi_i (s_i, s_{-i})$ for some $s_{-i} \in P_i[t_i]_{|S_{-i}}$.  But this contradicts that $(s_i,t_i)$ is rational. Thus, $R^{0}_{i|S_i} \subseteq S_i^{1}$.
Conversely, if $s_i \in S^{1}_i$, by rationality-completeness we have that, since $s_i$ is admissible with respect to $S_i \times S_{-i}$ there exists $t_i$ such that $\Pi_i(s_i, s_{-i}) \geq \Pi_{i}(\hat{s}_i, s_{-i})$ for any $s_{-i} \in P_{i}[t_i]_{|S_{-i}}$ and any $\hat{s}_i \in S_i$. Then $(s_i, t_i)$ is rational and belongs to $R^{0}_{i}$ and consequently, $s_i \in R^{0}_{i| S_i}$. That is, $R^{0}_{i|S_i} \supseteq S_i^{1}$.}
\item {\it Suppose the claim is true for $m = k$. Let us see that it is also true for $k+1$. 
Take $(s_i, t_i) \in R_i^{k+1}$. By definition, $(s_i, t_i)$ also belongs to $R_i^{k}$ and then $s_i \in S_i^{k+1}$. Suppose that $s_i \notin  S_i^{k+2}$. Then, there exists some $\hat{s}_i \in  S_i^{k+1}$ such that $\Pi_i (\hat{s}_i , s_{-i}) \geq \Pi_i (s_i, s_{-i})$ for some $s_{-i} \in S_{-i}^{k+2}$, but then, since $t_i \in R^{k+1}_{i | T_i}$, $s_i$ is a best response to any $s_{-i} \in R^{k}_{-i |S_{-i}}$$=$$S_{-i}^{k +1}$. In particular to any $s_{-i} \in S_{-i}^{k+2}$$\subseteq$$S_{-i}^{k+1}$. Contradiction. Then $ R_i^{k+1} \subseteq S_{i}^{k+2}$. 
For the converse, take $s_i \in S_i^{k+2}$ and notice that, by rationality-completeness, since $s_i$ is admissible with respect to $S_i^{k+1} \times S_{-i}^{k+1}$ there exists $t_i$ such that $\Pi_i(s_i, s_{-i}) \geq \Pi_{i}(\hat{s}_i, s_{-i})$ for any $s_{-i} \in P_{i}[t_i]_{|S^{k+1}_{-i}}$ and any $\hat{s}_i \in S^{k+1}_i$. Then $(s_i, t_i)$ belongs to $R^{k+1}_{i}$ and consequently, $s_i \in R^{k+1}_{i| S_i}$. That is, $R^{k+1}_{i|S_i} \supseteq S_i^{k+2}$}.
\end{itemize}

\section{Forcing in a Nutshell}

There are certain interesting properties that a particular $\mathbf{M}$ might exhibit. In particular, given the above result, we may be interested in the conditions that ensure that $\prod_{i}^n R^{\infty}_i$$\neq$$\emptyset$. Here is where the technique of {\em forcing} can be applied. 

To describe what {\em forcing} is we draw heavily from \cite{CG}, based in turn on the original work of Paul Cohen \cite{Co}. See also \cite{B}, \cite{Ba}, \cite{C},\cite{Je}, \cite{Je1} and \cite{K}, as well as the intriguing application in \cite{HW}.

Given a set-theoretic model $M$ of an axiomatically described system, we can define an extension of such a model by adjoining a set $\mathcal{G}$ to $M$, and denoting the new model by $M[\mathcal{G}]$. The nature of set $\mathcal{G}$, which we will call a {\it generic set} is such that, even being definable from within $M$, it is {\it indiscernable} from $M$. By this we mean that the language within $M$ allows us only to {\it name} the elements of $\mathcal{G}$, but not explicitly to describe their construction. In this way, we do not have access to the inner structure of $\mathcal{G}$, which remains unknowable from the point of view of $M$ -- hence the use of the word {\it generic} (as referring to the expression of something so ``mixed up'' or ``common'' that it cannot be discerned). Once the generic $\mathcal{G}$ has been defined, the {\it extension} via $\mathcal{G}$ of the ground model $M$ allows for new and possibly surprising ways to satisfy the ground axioms, with profound epistemic consequences. Indeed, although {\it truths} in $M[\mathcal{G}]$ are not directly accessible, we can define what we call a {\em forcing relation} between objects and relations at the level of the ground model. If one object {\em forces} a certain relation on $M$, then, if that relation belongs to $\mathcal{G}$ (and we might never be able to know that except as a modal claim across possible models), then we obtain ``truth'' in $M[\mathcal{G}]$. 

In what follows, we apply these methodological notions in order to extend strategic belief models as defined above, thereby making possible the forcing of a desirable epistemic property like $R\infty AR$ and indirectly the solution concept $\mathcal{S}_{IA}$ for {\em any} finite game $G$.

\section{Generic Types}
Let us start with a given $\mathbf{M}_0$ intended as a family of events of the game along with their underlying states of the game. This then constitutes our ground model, on which a $\mathbf{M}$ with desired properties will subsequently be built. We proceed by defining a family of forcing conditions, $\mathcal{P}$ with a partial order $\preceq$ defined as follows:

\begin{itemize}
\item $\mathcal{P}$$ = \{\phi =\langle (s_1, t_1), \ldots, (s_n, t_n)\rangle :$ every $i$ is rational and there exists a natural number $m(\phi)$ such that for any $i,j$, $t_j \in P^{m(\phi)}_{i}[t_i] \}$.
\item For any $\phi, \phi^{'} \in \mathcal{P}$, each one defined by a natural number ($m(\phi)$ and $m(\phi^{'})$), $\phi^{'}\preceq \phi$ iff $m(\phi) \geq m(\phi^{'})$.
\end{itemize}

We say that if $\phi^{'} \preceq \phi$, then $\phi$ {\em dominates} $\phi^{'}$. 

We then define the set of {\em correct} conditions, $\delta$, which satisfies the conditions of a {\em filter} in $\large( {\mathcal P}, \preceq \large)$:

\begin{itemize}
\item If $\phi^{'} \in \delta$ and $\phi^{'} \preceq \phi$ then $\phi \in \delta$. 
\item If $\phi^{'}, \phi^{''} \in \delta$ there exists $\phi \in \delta$ such that $\phi^{'} \preceq \phi$ and $\phi^{''} \preceq \phi$.
\end{itemize}

Our candidate is $\delta = \{\gamma = \langle (s_1, t_1), \ldots, (s_n, t_n)\rangle :$ there exists a natural number $m(\gamma)$ such that for any $i$, $(s_i, t_i) \in R^{m(\gamma)}_{i} \}$. We then have that:

\begin{prop}\label{correct}
$\delta$ is a correct set in $\mathcal{P}$.
\end{prop}
\noindent {\bf Proof}: {\it Let us see first that $\delta \subseteq \mathcal{P}$. That is, that for every $\gamma \in \delta$ there exists $\phi \in {\mathcal P}$ such that $\gamma = \phi$. We know that $\gamma$$=$$\langle (s_1, t_1), \ldots, (s_n, t_n)\rangle$ where for some $m \geq 0$ every $i$ is such that, $(s_i, t_i) \in R^{m}_{i}$. From this condition it follows that each $i$ is rational. We have now to show  that for every $i,j$, $t_j \in P^{m}_i[t_i]$. Suppose not. Then, there exists a $t_k$ such that $t_k \notin P^{m^{'}}_i[t_i]$, with $m^{'} \leq m$. This means that there exists $t_l \in P^{(m^{'}-1)}_i[t_i]$ such that $(s_k, t_k) \notin P_{l}[t_l]$. But, on the other hand, $(s_l, t_l) \in R^{m}_{l}$. We have by definition that $(s_l, t_l) \in S_l \times {\mathcal A}{\mathcal S}_{l}[R^{0}_k]$ which means that $t_l$ is such that $(s_k, t_k) \in P_l[t_l]$. Contradiction.} 

{\it The converse is also true: given a state $\langle (s_1, t_1), \ldots, (s_n, t_n)\rangle$$\in \mathcal{P}$, it follows that there exists $m$ such that each $(s_i, t_i) \in R^{m}_{i}$. Suppose not. Then, for a pair $i,j$, $t_j \notin {\mathcal A}{\mathcal S}_{j}[R^{m^{'}}_i]$ for some $m^{'} \leq m$. Then $P_j[t_j]_{|S_i \times T_i} \neq R^{m^{'}}_i$. In particular, we have that $t_i \notin P^{m^{'}}_j[t_j]$. Contradiction.}

{\it From this last implication it follows that if $\phi^{'} \in \delta$ and $\phi^{'} \preceq \phi$ then $\phi \in \delta$. This is because $m(\phi^{'}) \leq m(\phi)$ and $\phi$ is such that every $i$ is rational and for every pair $i,j$, $t_j \in P^{m(\phi)}_i$, which in turn implies that every $(s_i,t_i) \in R^{m(\phi)}_i$ and therefore, $\phi \in \delta$.}

{\it Finally, given $\phi^{'}, \phi^{''} \in \delta$, just take $m$ as the maximum of $m(\phi^{'})$ and $m(\phi^{''})$. Without loss of generality let us assume that $m = m(\phi^{'})$. Then, we take $\phi = \phi^{'}$ and it is easy to see that $\phi^{'} \preceq \phi$ and $\phi^{''} \preceq \phi$.}\\

We can now define a class of conditions called a {\em domination}, $D \subseteq {\mathcal P}$ (a {\em dense} set): 
\[\forall \phi^{'} \in {\mathcal P} \ \exists \phi \in D \ \mbox{such that} \ \phi^{'} \preceq \phi \]

We say that a correct set $\mathcal{G}$ is {\em generic} if $\mathcal{G} \subseteq \delta$ and $\mathcal{G} \cap D \neq \emptyset$ for any domination $D$. We then have that:

\begin{theo}
$\mathcal{G} = \{ \gamma \in \delta: \gamma =\langle (s_1, t_1), \ldots, (s_n, t_n)\rangle$ $\mbox{with for every} \ m $ $\mbox{and every} \ i,$ $(s_i, t_i) \in R^{m}_i \}$ is a generic set. 
\end{theo}
\noindent {\bf Proof}: {\it By Proposition~\ref{correct}, $\delta$ is a correct set and so is $\mathcal{G}\subset \delta$. To see that it is generic, just consider any $\phi \in \mathcal{P}$, which is identified by a finite natural number $m(\phi)$. Then, by definition, $\gamma \in \mathcal{G}$ is such that for every $i$, $(s_i, t_i) \in R^{m}_i$ for every $m$, in particular with $m \geq m(\phi)$. Then, $\phi \preceq \gamma$.}\\

In the context of the previous discussion, we can say that $\mathcal{G}$ defines a set of types $\{t^{*}_i\}_{i \in I}$ such that each one, joint with the corresponding $s^{*}_i$, satisfies the condition that $(s^{*}_{i}, t^{*}_{i}) \in R^{m}_{i}$ for {\em every} $m \geq 0$. That is, each $i$ is, with her type and the correponding strategy, rational and assumes rationality at all levels. In other words, it satifies $R\infty AR$. In accordance with Cohen's technique, these {\em generic types} cannot be defined in the language of $\mathbf{M_0}$. That is, there is no property $\lambda$ expressible in $\mathbf{M_0}$ such that:\footnote{The generic types are {\em indiscernable} in $\mathbf{M_0}$ \cite{CG}.} 
$$\forall t_i \in T_i \ \lambda(t_i) \Leftrightarrow \exists \gamma \in \mathcal{G} \ \mbox{such that} \ \exists s_i \in S_i \ (s_i, t_i) \in R^{m(\gamma)}_i$$ 

This realization is quite important since, as shown in \cite{BK}, $\mathbf{M_0}$ is not {\em definable complete}, i.e. there exists some event $E \in \mathbf{M_0}$, definable by a property $\lambda_E$  (i.e. $\langle (s_1, t_1), \ldots, (s_n, t_n) \rangle \in E \Leftrightarrow \lambda_{E}\large( \langle (s_1, t_1), \ldots, (s_n, t_n) \rangle \large)$) such that there exists a $i$ for whom no $t_i \in T_i$ satisfies $P_i[t_i] = E_{|\prod_{j \neq i} (S_j \times T_j)}$.\footnote{See also \cite{T05}.} 

But forcing shows that $\mathbf{M_0}$ can be extended to $\mathbf{M_0}[\mathcal{G}]$, to which the object constructed by $\mathcal{G}$, i.e. the class of generic types, belongs. To define $\mathbf{M_0}[\mathcal{G}]$, we will consider the {\em names} of objects in $\mathcal{G}$. The $\mathcal{G}$-names are recursively defined sets of the form $\{(\mu,\gamma):\mu \ \mbox{is a} \ \mathcal{G}-\mbox{name and} \ \gamma \in \mathcal{G}\}$. These can be ordered in terms of their rank, where a name $\mu$ of rank $0$ is the set of pairs $(\emptyset, \gamma)$ with $m(\gamma) =0$ and, recursively, we say that $\mu$ is of rank $m$ if it includes all the pairs $( \mu^{'}, \gamma)$ such that $m(\gamma) = m$ and the rank of $\mu^{'}$, $m^{'}$, verifies $m^{'} < m$.

The {\em referential value} of a name $\mu$, $r_{\mathcal{G}}(\mu)$ is then also defined recursively:

\begin{itemize}
\item If the rank of $\mu$ is $0$, $r_{\mathcal{G}}(\mu)= \{\langle (s_1, t_1), \ldots, (s_n, t_n) \rangle \in \prod_i (S_i \times T_i): (s_i, t_i) \in R^{0}_i \}$ iff there exists $(\emptyset, \gamma) \in \mu$. It is $r_{\mathcal{G}}(\mu) =  \mathbf{M_0}[\mathcal{G}]$, otherwise.
\item If the rank of $\mu$ is $m$, $r_{\mathcal{G}}(\mu)= \{ r_{\mathcal{G}}(\mu^{'}): \exists (\mu^{'}, \gamma) \in \mu\}$.
\end{itemize}

It is easy to see that names of rank $0$ yield all the states in which all players are rational, while for any $m > 0$ the names have as referential values all the states in which the players are rational and assume up to level $m$ the rationality of the others.

Then, $\mathbf{M_0}[\mathcal{G}]$ $=$ $\{r_{\mathcal{G}}(\mu): \mu \ \mbox{ is a name in} \ \mathbf{M_0}\}$. We have:

\begin{prop}
$\mathbf{M_0} \subset \mathbf{M_0}[\mathcal{G}]$.
\end{prop}
\noindent {\bf Proof}: {\it Trivial. Just take any name $\mu$ such that for every $\gamma \in \mathcal{G}$, $(\emptyset, \gamma) \notin \mu$. It is easy to find an event in $\mathbf{M_0}$ satisfying this condition: take anyone in which the states are such that there exist $i$ and $j$, with $(s_i, t_i) \in R^{m_i}_i$ and $(s_j, t_j) \in R^{m_j}_j$ and $m_i \neq m_j$.}\\

While $\mathcal{G}$ is not an element of $\mathbf{M_0}$, our recursive method has guaranteed that it exists in $\mathbf{M_0}[\mathcal{G}]$ as ``namable" although not strictly constructible.  This object represents the class of states in which each $(s_i, t_i; s_{-i}, t_{-i}) \in R\infty AR$. If we take the statement $\Psi(\langle (s_1, t_1), \ldots, (s_n, t_n) \rangle)$ which is true iff each $(s_i, t_i; s_{-i}, t_{-i}) \in R\infty AR$, we know that:

\[ {\mathbf M_0}[\mathcal{G}] \ \models \Psi(\langle (s_1, t_1), \ldots, (s_n, t_n) \rangle)\]

\noindent and by the central result in forcing\footnote{See \cite{B}, \cite{Je1}, \cite{C}, among others.}, if $\mu$ is a name in $\mathbf{M_0}$ such that $\langle (s_1, t_1), \ldots, (s_n, t_n) \rangle \in r_{\mathcal{G}}(\mu)$ we have that $\gamma \in \mathcal{G}$ is such that:

\[ \gamma \ \Vdash \Psi(\mu)\]

\noindent i.e. the generic types force $R\infty AR$ even though the space of types involved is not strictly complete.

We are now able to present the main result of the application of forcing to strategic belief models by bringing the discussion immediately above into conjunction with the earlier theorem concerning rationality-completeness:

\begin{prop}
If $\mathbf{M_0}$ is rationality-complete, $R^{\infty}_{i|S_i} = S_i^{\infty}$ for every $i$.
\end{prop}
\noindent {\bf Proof}: {\it Immediate from Theorem~\ref{admiss} and that ${\mathbf M_0}[\mathcal{G}]$ is such that $\prod_i^n R^{\infty}_{i}\neq \emptyset$.}\\

In this sense, we have defined a class of types involved in the characterization of $\mathcal{S}_{IA}$ over any $G$ with arbitrary (finite) numbers of players and actions. In other words, given the result above and our application of forcing to the space of types via generic types, we are able to ensure the sufficiency of rationality-completeness for the characterization of {\em IA} for a very general class of games.


\begin{thebibliography}{100}

\bibitem{AB95} Aumann, R. and Brandenburger, A.: Epistemic Conditions
for Nash Equilibrium, {\it Econometrica} {\bf 63}: 1161-1180, 1995.

\bibitem{B} Burgess, J.: Forcing, Chap.B4 in J. Barwise (ed.) {\bf Handbook of Mathematical Logic}, North-Holland, 1977.

\bibitem{Ba} Barwise, J. and Robinson, A.: Completing Theories by Forcing, {\it
Annals of Mathematical Logic} {\bf 2}: 119--142, 1970.

\bibitem{B} Brandenburger, A.: The Power of Paradox: Some Recent
Developments in Interactive Epistemology, {\it International Journal
of Game Theory} {\bf 35}: 465--492, 2007.

\bibitem{BF} Brandenburger, A. and Friedenberg, A.: Self-Admissible Sets, {\it Journal of Economic Theory} {\bf 145}: 785--811, 2009.

\bibitem{BK} Brandenburger, A. and Keisler, J.: An Impossibility
Theorem on Beliefs in Games, {\it Studia Logica} {\bf 84}: 211-240,
2006.

\bibitem{BFK} Brandenburger, A. - Friedenberg, A. - Keisler, J.: Admissibility in Games, {\em Econometrica} {\bf 76}: 307--352, 2008. 

\bibitem{C} Chow, T.: A Beginners Guide to Forcing, {\it Contemporary Mathematics} {\bf 479}: 25--40, 2009.

\bibitem{CG} Caterina, G. and Gangle, R.: Consequences of a Diagrammatic Representation of Paul Cohen's Forcing Technique Based on C.S. Peirce's Existential Graphs , {\it Studies in Computational Intelligence} {\bf 314}: 429--443, 2010.

\bibitem{Co} Cohen, P.: {\bf Set Theory and the Continuum Hypothesis}, Addison-Wesley, 1966.

\bibitem{HW} Hatchuel, A. and Weil, B.:  Design as Forcing: Deepening the Foundations of CK Theory, {\it Proceedings of the 16th International Conference on Engineering Design.} ({\tt http://www.cgs.ensmp.fr/design/}{\tt docs/Hatchuel-Weil-ICED07paperVdef21.pdf.}), 2007.

\bibitem{Je} Jech, T.: What is Forcing?, {\it Notices of the AMS} {\bf 55}: 692-693, 2008.

\bibitem{Je1} Jech, T.: {\bf Set Theory}, Springer Verlag, 2002.

\bibitem{K} Kanamori, A.: Cohen and Set Theory, {\it Bulletin of Symbolic Logic} {\bf 14}: 351--378, 2008.

\bibitem{KL} Keisler, J. and Lee, B.S.: Common Assumption of Rationality, {\it Working Paper}, 2011.

\bibitem{KM} Kohlberg, E. and Mertens, J-F.: On the Strategic Stability of Equilibria, {\it Econometrica} {\bf 54}: 1003--1037, 1986.

\bibitem{PR} Pacuit, E. and Roy, O.: Epistemic Foundations of Game Theory, {\it Working Paper}, 2012.

\bibitem{P} Perea, A.: {\bf Epistemic Game Theory: Reasoning and Choice}, Cambridge University Press, 2012.

\bibitem{T05} Tohm\'e, F.: Existence and Definability of States of the
World, {\it Mathematical Social Sciences} {\bf 49}: 81-100, 2005.



\end{thebibliography}
\end{document}